\documentclass[journal,twocolumn]{IEEEtran}
\usepackage[margin=1in]{geometry}
\usepackage{amsmath}
\usepackage{amsthm} 
\usepackage{amsfonts}
\usepackage{graphicx}
\graphicspath{{./images/}}
\newtheorem{lemma}{Lemma}
\newtheorem{theorem}{Theorem}

\newtheorem{definition}{Definition}
\newtheorem{remark}{Remark}

\newtheorem{assumption}{Assumption}

\usepackage{caption}
\usepackage{subcaption} 
\usepackage{xcolor}
\usepackage{algpseudocode}
\usepackage{algorithm}
\definecolor{kthblue}{rgb}{0.098, 0.329, 0.651}
\definecolor{newgreen}{RGB}{38, 115, 77}
\definecolor{newpurple}{RGB}{102, 0, 204}
\definecolor{newred}{RGB}{153, 0, 0}
\definecolor{newbrown}{RGB}{115, 77, 38}
\definecolor{neworange}{RGB}{255, 153, 0}
\usepackage{mathtools}
\usepackage[final]{changes}

\title{When are selector control strategies optimal for constrained monotone systems?}
\author{Hamed Taghavian, Ross Drummond and Mikael Johansson}

\begin{document}
\maketitle


\begin{abstract}
This paper considers optimal control problems defined by a monotone dynamical system, a monotone cost, and monotone constraints. We identify families of such problems for which the optimal solution is bang-ride, \emph{i.e.}, {that it always operates on the constraint boundaries and switches} between a finite number of state feedback controllers. This motivates the use of simpler policies, such as selector control, that can be designed without perfect models and full state measurements. The approach is successfully applied to several variations of the health-aware fast charging problem for lithium-ion batteries.
\end{abstract}



\section{INTRODUCTION}
Optimal control theory is a powerful framework for identifying the limits of performance {in dynamical systems and characterizing the optimal policies} that operate at these limits. {However, finding optimal feedback policies for constrained dynamical systems remains a major challenge in general.} It has become increasingly popular {to mimic these} policies using model-predictive controllers, but these are complex, memory-intensive, and require precise system models as well as significant online computations. For this reason, classical control designs remain highly popular in practice. In fact, it has been estimated that as much as $95\%$ of controllers deployed in the industry are PIDs~\cite{95pid}. Easy implementation and commissioning, simple tuning and analysis, and low requirements on the controller hardware are some of the advantages that allow classical controllers to remain competitive with more advanced options.

When applied to {constrained dynamical systems, classical controllers are often combined with a (min/max) selector, which is a nonlinear static device with multiple inputs and a single output~\cite{pidbook}.} {This device} uses simple logic to select one of the inputs, typically the one with the smallest or largest value, and connect it to the output. Hence, a selector essentially isolates individual control loops during operation, allowing for single-loop design methods to be applied in more complex control settings, thus simplifying the overall design.  We call a control system comprising a selector and multiple linear control loops a \emph{selector control} system.

Selector control systems are especially common in industrial applications where, for safety reasons, some variables must remain within certain limits, such as {the grid frequencies in power systems~\cite{pidbook}, the chemical concentrations in process control~\cite{Skog1}, and the desired thrust in engine control~\cite{engine1}.} Compared with more advanced methods, such as model predictive control, selector control is easier to implement and cheaper to execute. In addition, individual loops {in a selector control system} can often be tuned in a model-free fashion based on data. Despite these advantages, selector control is still applied in an ad-hoc manner, {with little theoretical understanding of the class of problems it performs well on~\cite{Skog1}}. Identifying {these} problems provides formal justification for {its} use in practice and can help improve {its} tuning.




{\subsection{Contributions}} 
{In this paper, we study optimal control problems} that have a monotone cost function, evolve with a monotone dynamical system, and are subject to monotone constraints. {We refer to this class as monotone optimal control problems.} While they may seem restrictive at first, monotone optimal control problems find practical applications in diverse areas such as molecular biology~\cite{mon_con_sys}, transportation networks~\cite{sad2016,ran2015}, rechargeable batteries~\cite{DCHCG2023}, and robotics~\cite{softrobots}. 

{We show that each monotone optimal control problem admits a unique bang-ride solution and provide sufficient conditions for this solution to be optimal (Theorem~\ref{thm:nonlin_decreasing}). Stronger results are obtained when the problem is linear. The necessary and sufficient condition for the optimality of bang-ride solutions in linear monotone optimal control problems is provided (Theorem~\ref{thm:NS_lin}). If in addition to linearity, the impulse response is decreasing, then we prove that the bang-ride policy is always optimal (Theorem~\ref{thm:Linear_bangride}). This latter condition can be verified from the step response alone without modeling the system. Since bang-ride solutions can be precisely mimicked by a selector control system, an important implication of our results is that} the optimal solutions to certain monotone control problems can be implemented using selector control on industrially relevant hardware. {We demonstrate our results by verifying the optimality of bang-ride protocols for charging lithium-ion battery cells.}





This paper is organized as follows: We introduce the family of monotone optimal control problems in Section~\ref{sec:monotone_optimal_control_problems} and discuss their solutions in {Sections~\ref{sec:solution_methods} and \ref{sec:BR}. An approximate implementation of the optimal solutions using selector control schemes is provided in Section~\ref{sec:selector_pid}.} A practical application of the results is provided in Section~\ref{sec:applications} and concluding remarks are given in Section~\ref{sec:conclusion}.

\subsection{Notation}
The $i^{\rm th}$ entry of the element $x_t\in\mathbb{R}^n$ in the sequence $x=\lbrace x_t\rbrace_{t=0}^{t=t_f}$ is denoted by $x_{i,t}$. All inequalities are understood componentwise, \emph{i.e.}, for two $n$-vectors $x_1$ and $x_2$, $x_1\geq x_2$ means that $x_{i,1}\geq x_{i,2}$ for all $i=1,2,\dots,n$.

\section{MONOTONE OPTIMAL CONTROL PROBLEMS}\label{sec:monotone_optimal_control_problems}

Consider the following discrete-time, possibly nonlinear, dynamical system with one input and multiple outputs 
\begin{align}\label{eqn:system}
\begin{array}{rcl}
    x_{t+1}&=&f(x_t,u_t) \\
    y_{t}&=&h(x_t,u_t) 
\end{array}& \quad  \mbox{ for } t\in\mathbb{N}_0.
\end{align}
Here, $x_t\in\mathbb{R}^n$ is the state, $y_{t}\in\mathbb{R}^p$ is the output vector and $u_t\in\mathbb{R}$ is the input. For a given initial state $x_0\in \mathbb{R}^n$, we study optimal control problems of the form
\begin{equation}
    \begin{array}[c]{rll}
    \underset{u}{\text{maximize}} & J(x_0,u)=\sum_{t=0}^{t_f}L(x_t,u_t) \\
    \mbox{subject to} & y_{t}\leq y_{\max},\\
    & u_t \leq u_{\rm max},\\
    & u, x, y \textnormal{ satisfy (\ref{eqn:system})}, \\
    \end{array} \label{eqn:optimization}
\end{equation}
for constant $u_{\max}\in \mathbb{R}$ and $y_{\max}\in \mathbb{R}^p$. For ease of notation, we have introduced $u=\lbrace u_t\rbrace_{t=0}^{t=t_f}$, 
$x=\lbrace x_t\rbrace_{t=0}^{t=t_f}$, and
$y=\lbrace y_t\rbrace_{t=0}^{t=t_f}$. When convenient, we also append $u_t$ to the output vector as
\begin{equation}\label{augmen}
    \bar{h}(x_t,u_t)=\begin{bmatrix}\bar{h}_0(x_t,u_t),\dots,\bar{h}_p(x_t,u_t)\end{bmatrix}^\top,
\end{equation}
where
    $$
    \bar{h}_i(x_t,u_t)=\left\lbrace\begin{array}{ll}
        u_t,  & i=0 \\
        h_i(x_t,u_t), & 1\leq i \leq p
    \end{array}\right.
    $$
and replace both the input and output inequality constraints in (\ref{eqn:optimization}) by the following equivalent constraint
$$
\bar{y}_{t}=\bar{h}(x_t,u_t)\leq \bar{y}_{\max}=\begin{bmatrix}u_{\max}& y^\top_{\max}\end{bmatrix}^\top.
$$

For a formal definition of monotone optimal control problems based on (\ref{eqn:optimization}), we need a few definitions and assumptions.

\begin{definition}
Let the $k$-step reachable set of (\ref{eqn:system}) from the initial state $z$ under the constraints of (\ref{eqn:optimization}) be denoted by
\begin{align*}
    \mathcal{R}(z,k) &=
    \lbrace x_k\in\mathbb{R}^n \,\vert\, x_0=z
     \mbox{ and } \nonumber\\
    &\quad\; x_t,u_t,y_{t}\textnormal{ satisfy (\ref{eqn:system})--(\ref{eqn:optimization})}\mbox{ for all }t\in[0,k] 
    \rbrace
\end{align*}
Then the cumulative reachable set of problem (\ref{eqn:optimization}) is defined as $\bar{\mathcal{R}}(x_0, t_f):= \bigcup_{k=0}^{t_f}  {\mathcal R}(x_0, k)$.
\end{definition}
The set $\bar{\mathcal{R}}(x_0, t_f)$ represents the set of states that can be visited on the time interval $[0, t_f]$ from initial state $x_0$ under the dynamics (\ref{eqn:system}) while satisfying the constraints of (\ref{eqn:optimization}).

\begin{definition}
Let ${\mathcal D}_x\subseteq \mathbb{R}^n$ and ${\mathcal D}_u\subseteq \mathbb{R}$. 
A function $g:\mathbb{R}^n\times \mathbb{R} \to \mathbb{R}^{n_g}$ is increasing in  $\mathcal{D}_x\times \mathcal{D}_u$ if for all $x_1,x_2\in\mathcal{D}_x$ and $u_1,u_2\in\mathcal{D}_u$ we have
$$
x_1\geq x_2\mbox{ and }u_1\geq u_2 \Rightarrow g(x_1,u_1)\geq g(x_2,u_2).   
$$
\end{definition}

A function $g$ is increasing if and only if all its components $g_1, \dots, g_{n_g}$ are increasing. We use the above definitions to impose monotonicity assumptions on the system dynamics (\ref{eqn:system}) and the optimization problem (\ref{eqn:optimization}).

\begin{assumption}[monotone cost]\label{ass:L}
The running-cost function $L(x_t,u_t)$ is differentiable and increasing in $(x_t,u_t)\in\mathbb{R}^n\times\mathbb{R}$.
\end{assumption}

\begin{assumption}[monotone system]\label{ass:f}
The transition function $f(x_t,u_t)$ is differentiable in $(x_t,u_t)\in\mathbb{R}^n\times\mathbb{R}$ and increasing in the region $(x_t,u_t)\in\mathcal{B}\times (-\infty,u_{\rm max}]$ where $\mathcal{B}\subset \mathbb{R}^n$ is a ball containing $\bar{\mathcal{R}}(x_0,t_f)$.
\end{assumption}

\begin{assumption}[monotone constraint]\label{ass:h}
The output function $h(x_t,u_t)$ is differentiable in $x_t\in\mathbb{R}^n$ for all $u_t\in(-\infty,u_{\rm max}]$, and both differentiable and strictly increasing in $u_t\in \mathbb{R}$ for all $x_t\in\bar{\mathcal{R}}(x_0,t_f)$.
\end{assumption}

We call Problem~(\ref{eqn:optimization}) a \emph{monotone optimal control problem} under Assumptions~\ref{ass:L}-\ref{ass:h}. We assume an optimal solution exists in this problem which may not be necessarily unique. The existence of an optimal solution in a feasible monotone optimal control problem can be ensured by, for example, Weierstrass theorem when the control inputs are constrained as $u_t\geq 0$ for all $t\in[0,t_f]$ (this condition is typically enforced in fast charging problems to avoid discharge).



Assumption~\ref{ass:h} states that 
$$
\partial h_i(x_t,u_t)/\partial u_t \geq 0, \quad i=1,2,\dots,p
$$
holds for all $x_t\in\bar{\mathcal{R}}(x_0,t_f)$ and $u_t\in \mathbb{R}$,
with equality met only on a set of zero measure. Therefore, since $\partial h_i(x_t,u_t)/\partial u_t \neq 0$ holds almost everywhere,
\begin{equation}\label{eqn:h1>h2}
    h(x_t,\hat{u}_t)>h(x_t,u_t)
\end{equation}
is satisfied whenever $\hat{u}_t>u_t$. Assumption~\ref{ass:h} is less restrictive than it seems because the inequality constraints associated with the output functions that are decreasing in $u_t$ can be safely removed from the optimization problem (\ref{eqn:optimization}) with no impact on the optimal value. To see why, we first note that by combining Assumptions~\ref{ass:L} and \ref{ass:f} one can show that $L(x_t,u_t)$ is increasing in $u_k$ for all $k\in\lbrace 0,1,\dots,t\rbrace$. Therefore, the cost function in (\ref{eqn:optimization}) is also increasing in $u_t$. Hence, there are optimal solutions to (\ref{eqn:optimization}) in which the inequality constraints associated with the output functions that are decreasing in $u_t$ are not activated (except in degenerate cases, \emph{e.g.}, where those constraints are explicitly enforced for all feasible $x$, such as the equality constraint $u_t=1$ enforced by the double inequality $h(x_t,u_t)=[-u_t,u_t]^\top \leq y_{\max}= [-1, 1]^\top$). When the output map does not depend on $u_t$, \emph{i.e.}, when $h_i(x_t,u_t)=h_i(x_t)$, 
Assumption~\ref{ass:h} is not satisfied. However, one may then redefine the $i$-th output to be
\begin{equation}\label{eqn:hhat}
    \hat{h}_i(x_t,u_t)=h_i\circ f^{d}_{ a}\circ f(x_t,u_t),
\end{equation}
where $f_a(\cdot)=f(\cdot,0)$ is the autonomous system dynamics in (\ref{eqn:system}), $f^{d}_a$ is the $d$-fold composition of $f_a$ with itself and
$d\in[0,n-1]\cup\lbrace \infty\rbrace$
is the relative degree of the $i$-th output in the nonlinear dynamic system (\ref{eqn:system}). Unless the output is completely decoupled from the control input, $d$ is finite, and hence, (\ref{eqn:hhat}) is well-defined. Recall that the relative degree $d$ determines the lag from the inputs to the outputs in a control system (\ref{eqn:system}). In particular, the first instance of the output signal which is affected by $u_0$ is $y_{d+1}=h_i(x_{d+1})$. Hence, by definition~(\ref{eqn:hhat}),
$$
\frac{\partial \hat{h}_i(x_t,u_t)}{\partial u_t} \not\equiv 0
$$
holds and the new output functions $\hat{h}_i$ satisfy Assumption~\ref{ass:h}. The transform in (\ref{eqn:hhat}) can be understood as a mere forward time shift
\begin{equation}\label{eqn:shifting_y}
    \hat{h}_i(x_t,u_t)=h_i(x_{t+d+1}).
\end{equation}
Therefore, by replacing the output functions $h_i$ by $\hat{h}_i$ in (\ref{eqn:optimization}), the same constraints are imposed on $y_{i,t}$ starting from $t=d+1$. {The initial steps prior to this time where the original output constraint is not enforced have a negligible impact when} the sampling time of the system is small (see \S~\ref{sec:applications}).

{\subsection{Related problems}}
The family of monotone dynamical systems defined by Assumptions~\ref{ass:f} and \ref{ass:h} is closely related to several other families of systems studied in the literature, such as incrementally positive systems~\cite{mon_con_sys}, monotone control systems~\cite{mon_con_sys}, and positive systems~\cite{1D2D}. In particular, the dynamical system (\ref{eqn:system}) is a monotone control system in the sense of~\cite{mon_con_sys} under Assumption~\ref{ass:f}, because 
$$
\hat{x}_0 \geq x_0\mbox{ and }\hat{u}_t\geq u_t\, \forall t\in\mathbb{N}_0
\Rightarrow
\hat{x}_t\geq x_t\, \forall t\in\mathbb{N}_0,
$$
which is readily shown by induction. If $f(0,0)=0$, $h(0,0)=0$ and  Assumption~\ref{ass:f} holds, the system (\ref{eqn:system}) is a discrete-time internally positive control system in the sense of~\cite{nonlin_pos}. {System (\ref{eqn:system}) can also describe all linear time-invariant internally positive control systems that have a single input and non-zero direct feedthroughs for each output.}

{
The optimal control problem~(\ref{eqn:optimization}) can also be viewed as a static optimization problem. To see this, one can simply remove the last constraint from (\ref{eqn:optimization}) by substituting the output functions from (\ref{eqn:system}), solving the state equations (\ref{eqn:system}) recursively for $t=0,1,\dots,t_f-1$, and finally, substituting the states $x_t$ in (\ref{eqn:optimization}). In this way, equation (\ref{eqn:optimization}) can be posed as the following static optimization problem in the variable $\vec{u}=[u_0,u_1,\dots,u_{t_f}]^\top\in\mathbb{R}^{t_f+1}$
\begin{equation}
    \begin{array}[c]{rll}
    \underset{\vec{u}}{\text{maximize}} & J(x_0,u)\\
    \mbox{subject to} & \vec{u}\in\mathcal{U}, \\
    \end{array} \label{eqn:static_optimization}
\end{equation}
It can be shown that when $\vec{u}\geq 0$ and $h(x_t,u_t)$ is an increasing function of $x_t\in\bar{\mathcal{R}}(x_0,t_f)$ for all $u_t\geq 0$, the objective function in (\ref{eqn:static_optimization}) is an increasing function of $\vec{u}$ by Assumptions~\ref{ass:L} and \ref{ass:f}, and the set $\mathcal{U}$ is normal by Assumptions~\ref{ass:f} and \ref{ass:h} (for the definition of a normal set, see~\cite{tuy1999normal}). In this case, problem (\ref{eqn:static_optimization}) can be posed as a monotonic optimization problem in the sense of~\cite{tuy2000monotonic}.}





\section{SOLUTION METHODS}\label{sec:solution_methods}

In this section, we study solutions to the optimization problem (\ref{eqn:optimization}) and propose a few different techniques to find and implement the optimal control input $u^{\star}$ in (\ref{eqn:optimization}).

\subsection{Open-loop sequences}\label{sec:open-loop}

The traditional way to approach discrete-time optimal control problems is to apply generic nonlinear programming techniques, such as sequential quadratic programming, to find a feasible sequence of control inputs $u=\lbrace u_t\rbrace_{t=0}^{t=t_f}$ that attains a large objective value. However, it is, in general, difficult to guarantee that an optimal solution will be found, or even to validate that a feasible control sequence is globally optimal. 

{For monotone optimal control problems (\ref{eqn:optimization}) in particular, it is also possible to apply the branch and bound procedures designed for monotonic optimization problems (\ref{eqn:static_optimization}), when the inputs are non-negative and the output functions are increasing in $x_t$~\cite{tuy2000monotonic}. The problem is simpler} when the transition function $f$ is linear. In this case, as long as the output functions $h_i(x)$ are convex and the running costs $L(x,u)$ are concave functions of their arguments, the optimal control problem (\ref{eqn:optimization}) can be posed as a standard convex optimization problem
with $(n+1)t_f+1$ variables and $(t_f+1)(p+1)+nt_f$ constraints. 
The optimal control sequence can then be computed efficiently using a wide range of reliable numerical routines.

{While the above approaches have some merits, they are only able to 
compute open-loop control sequences, making them sensitive to model uncertainties and hence too fragile to use in practice.} A practical example where this fragility could prove problematic is {fast charging lithium-ion batteries while accounting for cell degradation and model limitations (see Section~\ref{sec:applications}).} As the battery ages, its parameters and constraints will have to be adapted so as to ensure safety and performance, and enforcing this adaption efficiently requires a closed-loop approach. 


\subsection{Closed-loop policies}\label{sec:closed-loop}
To achieve more robust performance in practice, it is beneficial to {have a feedback policy where the control input is described in terms of the system states. One of the most popular forms of feedback policies used for monotone optimal control problems is the ``Bang-ride'' policy. This policy maximizes the control inputs to make at least one inequality constraint active in (\ref{eqn:optimization}) at \emph{every} time step.} {To express this policy as a feedback control law, we first note that the output functions are strictly increasing. Therefore each} {constraint equation
\begin{equation}\label{eqn:constraint_equation}
    \bar{h}_i(x_t,u_t)=\bar{y}_{i,\rm max}, \quad i=0,1,2,\dots,p
\end{equation}}
has at most one solution as $u_t=K_i(x_t)$, where $K_i:\mathbb{R}^n\to \mathbb{R}$ is a possibly nonlinear function. By defining $K_i(x_t)=+\infty$ in case equation (\ref{eqn:constraint_equation}) has no real solution for $u_t$, it is possible to write the bang-ride control law as the following nonlinear state-feedback
\begin{equation}\label{eqn:selector}
    u^{br}_t=\min\lbrace  K_0(x_t), K_1(x_t),
    \dots, K_p(x_t)\rbrace,
\end{equation}
where $K_0(x_t)\equiv u_{\rm max}$ corresponds to the only input constraint in (\ref{eqn:optimization}). As $K_0(x_t)$ is finite, the control law (\ref{eqn:selector}) is always well-defined. {Hence, every monotone optimal control problem (\ref{eqn:optimization}) admits one and only one bang-ride solution (\ref{eqn:selector}). The set of active inequality constraints defined by this control law is given by
\begin{equation}\label{eqn:I_t}
    \mathcal{I}^{br}_t:=\arg\min_{0\leq i\leq p}\lbrace K_i(x_t)\rbrace.
\end{equation}
Under the bang-ride operation, equality $\bar{y}_{i,t}=\bar{y}_{i,\rm max}$ holds for $i\in\mathcal{I}^{br}_t$.}

Implementing the bang-ride control law (\ref{eqn:selector}) is simple when an accurate model is available and all the system states
can be measured reliably. For example, in the linear case $\bar{h}(x_t,u_t)=\bar{C}x_t+\bar{D}u_t$ where $\bar{C}\in\mathbb{R}^{(p+1)\times n}$ and $\bar{D}\in\mathbb{R}^{p+1}$ are constant, the feedback functions $K_i$'s are simply given by
\begin{align*}
K_i(x_t)&=(\bar{y}_{i,\rm max}-\bar{C}_ix_t)/\bar{D}_i, \quad i=0,2,\dots,p
\end{align*}
where $\bar{C}_i$ and $\bar{D}_i$ denote the $i$th rows of matrices $\bar{C}$ and $\bar{D}$. In the general case, if it is difficult to find $K_i$ analytically, the monotonicity of the output functions facilitates solving (\ref{eqn:constraint_equation})
numerically for a given $x_t$ to determine each $K_i(x_t)$. Unlike the analytic approach, which can be executed offline, the numerical approach requires that the equations (\ref{eqn:constraint_equation}) are solved online in every time step.




Although the optimality of bang-ride policies has been proven for specific systems~\cite{PLATM2020} and conjectured for others~\cite{LVWWZTWZ2021}{, it is not clear whether these policies are optimal in general. Therefore, we investigate the optimality of bang-ride solutions in the next section.}

\section{OPTIMALITY OF BANG-RIDE POLICIES}\label{sec:BR}


{We start this section with a simple counter-example which shows that bang-ride solutions are \emph{not} always optimal in monotone optimal control problems.} Consider an instance of (\ref{eqn:optimization}) where the dynamics and running costs are linear, \emph{i.e.},
\begin{align}\label{eqn:lin}
    f(x_t,u_t)&=Ax_t+Bu_t \nonumber\\
    h(x_t,u_t)&=Cx_t+Du_t \\
    L(x_t,u_t)&=Ex_t+F u_t\nonumber
\end{align}
 and the matrices in (\ref{eqn:lin}) are all non-negative. Furthermore, let $x_0=0$, $t_f=1$, $u_{\rm max}>0$, $y_{\rm max}=D u_{\rm max}$, $E=0$ and $F=1$. The constraints for this problem are
\begin{align*}
  &\left( u_0\leq u_{\rm max} \right)\,\wedge\,\left(y_0=D u_0\leq y_{\rm max}\right)  &:t=0\\
  &\left( u_1\leq u_{\rm max} \right)\,\wedge\,\left(y_1=CB u_0 + D u_1\leq y_{\rm max}\right)&: t=1
\end{align*}
The two constraints at $t=0$ are equivalent, and a bang-ride solution $u^{br}=\lbrace u_0,u_1\rbrace$ must have $u_0=u_{\rm max}$. We can therefore re-write the constraints at $t=1$ as
\begin{equation}\label{eqn:u1ex}
    \left(u_1\leq u_{\rm max} \right)\,\wedge\,\left(D u_1\leq (D-CB) u_{\rm max} \right).
\end{equation}
Define
$$
\gamma=\min_i \lbrace (D_i-C_iB)/D_i \rbrace,
$$
where $C_i$ and $D_i$ denote the $i$th rows of matrices $C$ and $D$. Since the system is positive, we have $\gamma \leq 1$. Therefore, $u_1=\gamma u_{\rm max}$ holds for any bang-ride solution satisfying (\ref{eqn:u1ex}). The total cost of $u^{br}$ is
$$
J(0,u^{br})=u_0+u_1=(1+\gamma)u_{\rm max}.
$$
Now consider another control sequence  $u=\lbrace 0,u_{\rm max}\rbrace$. It is easy to show that this control sequence is feasible and $J(0,u)=u_{\rm max}$. Comparing the two costs yields
$$
J(0,u^{br})-J(0,u)=\gamma u_{\rm max}.
$$
Hence, when $\gamma<0$, the bang-ride control is not optimal.

{Since the bang-ride solution (\ref{eqn:selector}) is not always optimal, the inequality constraints in (\ref{eqn:optimization}) may not be activated at every time instant during the optimal operation. When the inputs are non-negative and the output functions are increasing in $x_t$, the problem (\ref{eqn:static_optimization}) becomes a monotonic optimization problem which has an optimal solution on the boundary of the constraint set $\mathcal{U}$ (see Proposition~7 in \cite{tuy2000monotonic}). This proves an inequality constraint in (\ref{eqn:optimization}) is active at \emph{some} time instant during the optimal operation. However, it is not clear at which time instants the inequality constraints become active and whether this result is true in the general case when the output functions are not increasing. Therefore, we propose the following lemma which shows an inequality constraint is always active at the \emph{final} time instant during the optimal operation. This result is essential to obtain the conditions that ensure the optimality of the bang-ride solutions in Theorem~\ref{thm:nonlin_decreasing}.} 


\begin{lemma}\label{thm:bangride} 
Under Assumptions~\ref{ass:L}--\ref{ass:h}, the following two statements hold for the optimal control problem \eqref{eqn:system}-\eqref{eqn:optimization}: 
\begin{itemize}
    \item [1)] There is an optimal solution $u^{\star}$ to problem~(\ref{eqn:optimization}) that cannot be increased at any single time $t\in [0, t_f]$.
    \item [2)] Under this optimal solution, at least one inequality constraint in (\ref{eqn:optimization}) is active at $t=t_f$. 
\end{itemize}

\end{lemma}
\begin{proof}
Let $u^{\star}$ be an optimal solution to the optimal control problem  \eqref{eqn:system}-(\ref{eqn:optimization}). 
Define another input sequence $\hat{u}(\varepsilon)=\lbrace \hat{u}_t\rbrace_{t=0}^{t=t_f}$ by increasing $\hat{u}^{\star}$ at a single time $\tau\in [0, t_f]$ as
\begin{equation}\label{eqn:uhat}
    \hat{u}_t=\left\lbrace
    \begin{array}{ll}
    u^{\star}_t + \varepsilon,  &  t=\tau \\
    u^{\star}_t,  & t\neq \tau
    \end{array}\right.,
\end{equation}
where $\varepsilon \geq 0$. Let $x^{\star}=\lbrace x^{\star}_t\rbrace_{t=0}^{t=t_f}$ and $\hat{x}(\varepsilon)=\lbrace \hat{x}_t\rbrace_{t=0}^{t=t_f}$ denote the state trajectories of the system (\ref{eqn:system}) associated with the control inputs $u^{\star}$ and (\ref{eqn:uhat}) respectively. As system (\ref{eqn:system}) is causal, the trajectory $\hat{x}(\varepsilon)$ associated with $\hat{u}(\varepsilon)$ satisfies
\begin{align}
\hat{x}_t &=x^{\star}_t, \quad t=0,1,\dots,\tau. \nonumber
\intertext{We will now use induction to show that the monotonicity of the system (Assumption~\ref{ass:f}) implies that}
    \hat{x}_t &\geq x^{\star}_t, \quad t=\tau+1, \dots, t_f. \label{eqn:xhat>xstar}
\end{align}
To this end, first note that since $\hat{u}_{\tau}\geq u^{\star}_{\tau}$ and the function $f$ is increasing in the second argument, inequality (\ref{eqn:xhat>xstar}) holds for $t=\tau+1$ because
$$
\hat{x}_{\tau+1}=f(\hat{x}_{\tau},\hat{u}_{\tau})\geq f(x^{\star}_{\tau},u^{\star}_{\tau})=x^{\star}_{\tau+1}.
$$
Moreover, if (\ref{eqn:xhat>xstar}) holds for some $t\in[\tau+1,t_f-1]$, then 
$$
\hat{x}_{t+1}=f(\hat{x}_t,\hat{u}_t)\geq f(x^{\star}_t,u^{\star}_t)=x^{\star}_{t+1},
$$
since $f$ is increasing in its first argument. This proves that (\ref{eqn:xhat>xstar}) holds for all $\tau+1\leq t \leq t_f$. 
Next, define
\begin{align}\label{eqn:J>=J}
        \Delta(\varepsilon)&=J(x_0,\hat{u}(\varepsilon))-J(x_0,u^{\star})\nonumber\\
        &=\sum_{t=0}^{t_f} L(\hat{x}_t,\hat{u}_t)-\sum_{t=0}^{t_f} L(x^{\star}_t,u^{\star}_t). 
\end{align}
The monotonicity of running costs implies $\Delta(\varepsilon)\geq~0$, while the optimality of $u^{\star}$ gives  $\Delta(\varepsilon)\leq 0$. Therefore, we have $\Delta(\varepsilon)= 0$. As this equality holds for $\varepsilon=0$ (because $\hat{u}(0)=u^{\star}$) and that $\Delta(\varepsilon)$ is continuous {and increasing} in $\varepsilon$ (by Assumptions~\ref{ass:L} and \ref{ass:f}), there exists $\varepsilon_{\max}\in [0, u_{\max}-u^{\star}_{\tau}]$ as the largest value of $\varepsilon$ 
such that $\hat{u}(\varepsilon)$ is feasible and  $\Delta(\varepsilon)=0$ for all $\varepsilon \in [0,\varepsilon_{\max}]$. Clearly, $\hat{u}(\varepsilon_{\max})$ is also optimal and cannot be increased at time $\tau$ without violating the constraints. By repeating the above process with $u^\star=\hat{u}(\varepsilon_{\max})$
as the starting point, one can find an optimal solution that is maximal at every time step. This proves point~1.

%
%

Now let $u^{\star}$ be the optimal solution that is maximal at every time step. For any $\tau\in\lbrace 0,1,\dots,t_f\rbrace$, an increment in $u^{\star}_{\tau}$ (while keeping the equality constraints (\ref{eqn:system}) satisfied and keeping $u^{\star}_t$ constant for all $t\neq \tau$) will violate at least one constraint in the optimal control problem (\ref{eqn:optimization}). Since the system (\ref{eqn:system}) is causal, the violation  must occur at some $t^{\star}\in\lbrace \tau,\tau+1,\dots,t_f\rbrace$. Since the function $f$ is continuous (Assumption~\ref{ass:f}), the increment in $u^{\star}_{\tau}$ and $x_{t^{\star}}^{\star}$ can be chosen arbitrarily small. As the output function $h$ is also continuous (Assumption~\ref{ass:h}), the increments in $y_{i,t^{\star}}=h_i(x^{\star}_{t^{\star}},u^{\star}_{t^{\star}})$ can be made as small as desired, while violating one of the inequality constraints
in (\ref{eqn:optimization}). This indicates that the corresponding constraint is active under the optimal control sequence $u^{\star}$. The active constraint time $t^\star$ can be uniquely identified for $\tau=t_f$ as $t^{\star}=t_f$. This proves point~2. 
 

\end{proof}

For certain families of monotone optimal control problems, an optimal solution not only activates an inequality constraint at the \emph{final} time instant (Lemma~\ref{thm:bangride}) but it also activates at least one inequality constraint at \emph{every} time instant $t\in [0, t_f]$. {The optimal solution to these problems is given by a feedback policy described by the bang-ride control law (\ref{eqn:selector}). The rest of this section is devoted to identifying these problems.}

\subsection{Nonlinear systems}
{As discussed earlier, in addition to Assumptions~\ref{ass:L}-- \ref{ass:h}, the monotone optimal control problem
(\ref{eqn:optimization}) needs to satisfy some extra conditions to ensure the bang-ride feedback policy (\ref{eqn:selector}) is optimal. The following theorem gives such conditions for nonlinear systems based on Lemma~\ref{thm:bangride}.}

\begin{theorem}\label{thm:nonlin_decreasing}
Let Assumptions~\ref{ass:L}--\ref{ass:h} hold, the output function $h(x_t,u_t)$ is decreasing in $x_t$, and $\partial h_i(x_t,u_t)/\partial u_t\neq 0$, where $(x_t,u_t)\in\bar{\mathcal{R}}(x_0,t_f)\times(-\infty,u_{\rm max}]$. {Then the bang-ride solution (\ref{eqn:selector}) is optimal in the monotone optimal control problem (\ref{eqn:optimization}).}
\end{theorem}
\begin{proof}
    Let $\mathcal{I}_t$ denote the set of indices $i$ such that {$\bar{y}_{i,t}=\bar{y}_{i,\rm max}$} holds. We will use induction to prove that {there is an optimal solution $u^{\star}$ of problem~(\ref{eqn:optimization}) that makes at least one inequality constraint active at every time step $t\in[0,t_f]$, that is,
    \begin{equation}\label{eqn:Inotempty}
            \mathcal{I}_t\neq \emptyset,\quad \forall t\in [0, t_f].
    \end{equation}    
    This will prove the bang-ride control law (\ref{eqn:selector}) is optimal in (\ref{eqn:optimization}), as it is the only control sequence satisfying the property (\ref{eqn:Inotempty}).} 
    
    First, we may assume that there is at least one inequality constraint active at the final time under the optimal sequence $u^{\star}$ according to Lemma~\ref{thm:bangride}, \emph{i.e.},
$\mathcal{I}_{t_f}\neq \emptyset$. Now assume
    \begin{equation}\label{eqn:Itnotempty_nonlin}
        \mathcal{I}_t\neq \emptyset,\quad \forall t\in[\tau+1,t_f]
    \end{equation}
    for some $\tau \in [0,t_f-1]$. We will show that assumption (\ref{eqn:Itnotempty_nonlin}) implies there is also an inequality constraint active at time $t=\tau$, \emph{i.e.},
    $$
    \mathcal{I}_{\tau}\neq \emptyset.
    $$
    Define a control sequence $\hat{u}(\varepsilon)=\lbrace\hat{u}_t\rbrace_{t=0}^{t=t_f}$ by perturbing the optimal solution $u^{\star}$ as
    \begin{equation}\label{eqn:uhat_nonlin_decreasing}
    \hat{u}_t=\left\lbrace\begin{array}{ll}
    u^{\star}_t+\varepsilon_t,& t\geq \tau \\
    u^{\star}_t,& t< \tau 
    \end{array}\right.,
    \end{equation}
    where $\varepsilon_{{\tau}}>0$. Let the perturbations $\varepsilon\in\mathbb{R}^{t_f-\tau+1}$ be chosen such that the active constraints in the range $t\in[\tau+1,t_f]$ remain active under the new control sequence $\hat{u}(\varepsilon)$. This happens if for all $t\in[\tau+1,t_f]$ and $i\in\mathcal{I}_t$,
    \begin{align}\label{eqn:dy=0}
    &d\bar{y}_{i,t}=\\
    &\frac{\partial \bar{h}_i(x_t,u_t)}{\partial u_t}\varepsilon_t+
    \frac{\partial \bar{h}_i(x_t,u_t)}
    {\partial x_t}\frac{\partial f(x_{t-1},u_{t-1})}
    {\partial u_{t-1}}\varepsilon_{t-1}+\nonumber\\
    &\sum_{k=\tau}^{t-2} \frac{\partial \bar{h}_i(x_t,u_t)}
    {\partial x_t} \left(\prod_{j=t-1}^{k+1}\frac{\partial f(x_{j},u_{j})}{\partial x_{j}}\right)
    \frac{\partial f(x_{k},u_{k})}{\partial u_{k}}\varepsilon_{k}\nonumber\\
    &=0.\nonumber
    \end{align}
    as $\varepsilon\to 0$, where all the partial derivatives are evaluated along the optimal trajectory $x=x^{\star}$ and $u=u^{\star}$. We show that equation (\ref{eqn:dy=0}) yields
    \begin{equation}\label{eqn:eps>0}
            \varepsilon_k \geq 0, \quad k=\tau,\tau+1,\dots,t_f
    \end{equation}   
    by induction. Starting with $t=\tau$, one has $\varepsilon_k>0$ by definition (\ref{eqn:uhat_nonlin_decreasing}). Assuming that (\ref{eqn:eps>0}) holds for all $k\in\lbrace\tau,\dots,t-1\rbrace$, one has
    \begin{align*}
    &\frac{\partial \bar{h}_i(x_t,u_t)}{\partial u_t}\varepsilon_t=\\
    &-\sum_{k=\tau}^{t-1} \frac{\partial \bar{h}_i(x_t,u_t)}
    {\partial x_t} \left(\prod_{j=t-1}^{k+1}\frac{\partial f(x_{j},u_{j})}{\partial x_{j}}\right)
    \frac{\partial f(x_{k},u_{k})}{\partial u_{k}}\varepsilon_{k}
    \end{align*}
    from (\ref{eqn:dy=0}), which is non-negative due to Assumption~\ref{ass:f} and the fact that {$\bar{h}(x_t,u_t)$} is decreasing in $x_t$. Hence, as $\partial \bar{h}_i(x_t,u_t)/\partial u_t>0$, one obtains $\varepsilon_t\geq 0$. This proves (\ref{eqn:eps>0}). Now consider the first-order perturbation of the cost function as
    \begin{align}\label{eqn:Cost_purt}
    \Delta(\epsilon)&=J(x_0,\hat{u}(\varepsilon))-J(x_0,u^{\star})=\nonumber\\
    &=\sum_{t=\tau}^{t_f} \frac{\partial L(x_t,u_t)}{\partial u_t}\varepsilon_t
    \nonumber\\
    &+\sum_{k=0}^{t_f-\tau-1}\sum_{t=\tau}^{t_f-k-1}
    \frac{\partial L(x_{t+k+1},u_{t+k+1})}{\partial x_{t+k+1}}\times\nonumber\\
    &\left(\prod_{j=t+k}^{j=t+1}\frac{\partial f(x_{j},u_{j})}{\partial x_{j}}\right)
    \frac{\partial f(x_{t},u_{t})}{\partial u_{t}}\varepsilon_{t}.
    \end{align}
    where all the partial derivatives are evaluated along the optimal trajectory $x=x^{\star}$ and $u=u^{\star}$. Since all the terms in (\ref{eqn:Cost_purt}) are non-negative by (\ref{eqn:eps>0}) and Assumptions~\ref{ass:L} and \ref{ass:f},
    \begin{equation}\label{eqn:Deltalin>=0_first}
    \Delta(\varepsilon)\geq 0
    \end{equation}
    holds true. Two cases are possible based on (\ref{eqn:Deltalin>=0_first}). \\  
    \noindent\textit{Case 1}: $J(x_0,\hat{u}(\varepsilon))>J(x_0,u^{\star})$ holds in a deleted neighborhood of $\varepsilon=0$. In this case $\hat{u}(\varepsilon)$ is infeasible (since $u^{\star}$ is assumed to be an optimal solution), which means $\hat{u}(\varepsilon)$ violates at least one inequality constraint. Since system (\ref{eqn:system}) is causal, the violated constraint must occur
    at some $t^{\star}\in\lbrace\tau,\dots,t_f\rbrace$. However, as the active constraints in the range $t\in\lbrace\tau+1,\dots,t_f\rbrace$ are kept satisfied under the perturbed input sequence $\hat{u}(\varepsilon)$, the violated constraint must occur at $t=\tau$. {Therefore, an inequality constraint is active at $t=\tau$ under $u^\star$.}\\
    \noindent\textit{Case 2}: $J(x_0,\hat{u}(\varepsilon))=J(x_0,u^{\star})$ holds in a neighborhood of $\varepsilon=0$. In this case, the sequence $\hat{u}(\varepsilon)$ is also an optimal solution to Problem~\ref{eqn:optimization}, which satisfies
    $$
    \hat{u}_{\tau}>u^{\star}_{\tau}.
    $$
    Using the new optimal sequence $\hat{u}(\varepsilon)$ as the starting point and repeating the same process described above, one can obtain an optimal sequence that makes an inequality constraint active at $t=\tau$, because $h_i(x_{\tau},u_{\tau})$ is strictly increasing in $u_\tau$ for all $i$ (Assumption~\ref{ass:h}).
    
    {We have shown that if there is an inequality constraint active at every $t\in[\tau+1,t_f]$ under an optimal operation, an inequality constraint can also be activated at $t=\tau$ (if it is not already active) without losing optimality. Hence, there is an optimal solution satisfying (\ref{eqn:Inotempty}) and the bang-ride control law is optimal.
    }
\end{proof}

\subsection{Linear systems}

{Theorem~\ref{thm:nonlin_decreasing} only provides a sufficient condition for optimality of bang-ride solutions. However, when the dynamic system and the running-cost function are linear as (\ref{eqn:lin}), one can obtain the necessary and sufficient condition of optimality for the bang-ride solution (\ref{eqn:selector}). For a convenient expression of this condition, we define the auxiliary system
\begin{align}\label{eqn:aux}
\begin{array}{rcl}
    x_{t+1}&=&Ax_t+Bu_t \\
    z_{t}&=&Ex_t+Fu_t 
\end{array}& \quad  \mbox{ for } t\in\mathbb{N}_0,
\end{align}
by merging the running cost and the state dynamics in (\ref{eqn:lin}). We let $s_t$ be the step response of the auxiliary system (\ref{eqn:aux}) and
\begin{equation}\label{eqn:gt}
\bar{g}_t=\left\lbrace\begin{array}{ll}
     \bar{D},&  t=0\\
     \bar{C} A^{t-1}B,& t\geq 1 
\end{array}\right.
\end{equation}
be the impulse response of the original system (\ref{eqn:system}) with the augmented output (\ref{augmen}), where
$$
\bar{C}=\begin{bmatrix}
    0 \\ C
\end{bmatrix}, \quad
\bar{D}=\begin{bmatrix}
    1 \\ D
\end{bmatrix}.
$$
\begin{theorem}\label{thm:NS_lin}
    Let Assumptions~\ref{ass:L}-- \ref{ass:h} hold, the functions $f$, $h$, $L$ be linear as (\ref{eqn:lin}), and $\mathcal{I}^{br}_t=\lbrace i^{br}_t\rbrace$ be singleton for all $t\in[0,t_f]$. Define the matrices $\bar{G}\in\mathbb{R}^{(t_f+1)\times(t_f+1)}$ and $S\in\mathbb{R}^{1\times (t_f+1)}$ as
    $$
    [\bar{G}]_{ij}=\left\lbrace 
    \begin{array}{ll}
    \bar{g}_{i^{br}_{i-1},i-j},& j\leq i \\
    0, & j>i
    \end{array}
    \right.
    ,\quad [S]_{1,i}=s_{t_f-i+1}
    $$
    for $i,j\in\lbrace 1,2,\dots,t_f+1\rbrace$. Then the bang-ride control law (\ref{eqn:selector}) is an optimal solution of (\ref{eqn:optimization}) if and only if
    \begin{equation}\label{eqn:SG>0}
        S\bar{G}^{-1}\geq 0.
    \end{equation}
\end{theorem}
\begin{proof}
    Consider perturbing the bang-ride control sequence $u^{br}$ as
    $$
    \hat{u}_t=u^{br}_t+\varepsilon_t,\quad t=0,1,\dots,t_f.
    $$
    The perturbed control sequence $\hat{u}(\varepsilon)=\lbrace\hat{u}_t\rbrace_{t=0}^{t=t_f}$ is feasible in (\ref{eqn:optimization}) as $\varepsilon\to 0$, if and only if
    \begin{align}\label{eqn:dy=0_lin}
    d \bar{y}_{i,t}&=\bar{D}_i\varepsilon_t
    +\sum_{k=0}^{t-1} \bar{C}_iA^{t-k-1}B \varepsilon_{k}\nonumber\\
    &=\sum_{k=0}^{t} \bar{g}_{i,t-k}\varepsilon_k\leq 0 
    \end{align}
    holds true for all $t\in[0,t_f]$ and $i\in\mathcal{I}^{br}_t$. This is the case, if there is some $\zeta\in\mathbb{R}_{\geq 0}^{t_f+1}$ such that
    \begin{equation}\label{eqn:eps=Gzeta}
        \varepsilon=-\bar{G}^{-1}\zeta,
    \end{equation}
    where $\varepsilon=\begin{bmatrix}\varepsilon_0,\dots,\varepsilon_{t_f}\end{bmatrix}^\top$. The matrix $\bar{G}$ is invertible as
    $$
    \det(\bar{G})=\prod_{t=0}^{t_f}[\bar{D}]_{i^{br}_t}>0
    $$
    holds under Assumption~\ref{ass:h}. The differential of the cost function with respect to the feasible changes (\ref{eqn:eps=Gzeta}) in the bang-ride control input is given by
    \begin{align}\label{eqn:Jsigma} 
    dJ\bigl(x_0,u^{br}\bigr)&=\sum_{t=0}^{t_f} F\varepsilon_t+\sum_{t=0}^{t_f}\sum_{k=0}^{t-1}EA^{t-1-k}B\varepsilon_{k} \nonumber\\
    &=\sum_{k=0}^{t_f} s_{t_f-k}\varepsilon_{k},
    \end{align}
    which, after substituting (\ref{eqn:eps=Gzeta}), results in
    \begin{equation}
        dJ\bigl(x_0,u^{br}\bigr)=-S\bar{G}^{-1}\zeta.
    \end{equation}
    The bang-ride sequence $u^{br}$ is a local maximum in the optimization problem (\ref{eqn:optimization}) if and only if $dJ\bigl(x_0,u^{br}\bigr)\leq 0$ holds in all feasible directions, that is, for all $\zeta\in\mathbb{R}_{\geq 0}^{t_f+1}$. This is true if and only if the non-negative orthant is invariant under $S\bar{G}^{-1}$, which is equivalent to condition (\ref{eqn:SG>0}). Furthermore, under the condition (\ref{eqn:lin}), problem (\ref{eqn:optimization}) is a linear program. Therefore, the local maximality of $u^{br}$ implies that (\ref{eqn:selector}) is also the (global) optimal solution of (\ref{eqn:optimization}).
\end{proof}}

\begin{remark}{The condition that $\mathcal{I}^{bt}_t$ is singleton ($\vert\mathcal{I}^{bt}_t\vert=1$) in Theorem~\ref{thm:NS_lin} is not much restrictive as the feedback control functions $K_i(x_t)$ in (\ref{eqn:selector}) are only compared numerically up to machine precision. Nevertheless, this condition can be dropped by generalizing Theorem~\ref{thm:NS_lin} to include the problems where multiple inequality constraints become active at the same time in (\ref{eqn:optimization}) during the bang-ride operation, \emph{i.e.},
$$
\vert\mathcal{I}^{bt}_t\vert\geq 1.
$$
In this case, one can redefine $\bar{G}\in\mathbb{R}^{(\sum_{t=0}^{t_f}\vert\mathcal{I}^{br}_t\vert)\times(t_f+1)}$ as follows
    $$
    \bar{G}=
    \begin{bmatrix}
        \bar{g}_{\mathcal{I}^{br}_0,0}& 0 & \cdots & 0\\
        \bar{g}_{\mathcal{I}^{br}_1,1}& \bar{g}_{\mathcal{I}^{br}_1,0} & \cdots & 0\\
        \vdots & \ddots & &\vdots \\
        \bar{g}_{\mathcal{I}^{br}_{t_f},t_f}& \bar{g}_{\mathcal{I}^{br}_{t_f},t_f-1} & \cdots & \bar{g}_{\mathcal{I}^{br}_{t_f},0}
    \end{bmatrix},
    $$
    where
    $$
    \bar{g}_{\mathcal{I}^{br}_t,\tau}=\begin{bmatrix}  \bar{g}_{i_1,\tau},\bar{g}_{i_2,\tau},\dots,\bar{g}_{i_{\vert\mathcal{I}^{br}_t\vert},\tau} 
    \end{bmatrix}^\top
    $$
    and $\mathcal{I}^{br}_t=\lbrace i_1,\dots,i_{\vert\mathcal{I}^{br}_t\vert}\rbrace$. Let us denote the intersection of the image of $\bar{G}$ and the non-negative orthant as the polyhedral convex cone
    $$\mathcal{K}=\left\lbrace \left.Z\alpha\,\right\vert\, \alpha\in\mathbb{R}^k_{\geq 0}\right\rbrace.
    $$
    Then the necessary and sufficient condition (\ref{eqn:SG>0}) is replaced by
    $$
    S\bar{G}^{\dagger}Z\geq 0,
    $$
    where $\bar{G}^{\dagger}=(\bar{G}^\top\bar{G})^{-1}\bar{G}^\top$.
}\end{remark}


{
In order to check the condition provided in Theorem~\ref{thm:NS_lin}, one needs to have access to the full model and simulate the step response of an auxiliary system (\ref{eqn:aux}). Instead, in the next theorem, we provide a condition that can be verified by just observing the step response of the \emph{original} system (\ref{eqn:system}). This condition is however only sufficient.
}

\begin{theorem}\label{thm:Linear_bangride}
{Let Assumptions~\ref{ass:L}-- \ref{ass:h} hold and the functions $f$, $h$, $L$ be linear as (\ref{eqn:lin}). If the step response of (\ref{eqn:system}) is \emph{concave} (\emph{i.e.}, the impulse response is decreasing) in the range $t\in[0,t_f]$, 
then the bang-ride control law (\ref{eqn:selector}) is optimal in the monotone optimal control problem (\ref{eqn:optimization}).}
\end{theorem}
\begin{proof}
    The proof follows a similar set of arguments as the proof of Theorem~\ref{thm:nonlin_decreasing}. Let $\mathcal{I}_t$ denote the set of indices $i$ such that $\bar{y}_{i,t}=\bar{y}_{i,\rm max}$ holds. We will use induction to prove that there is an optimal solution $u^{\star}$ of problem~(\ref{eqn:optimization}) which makes $\mathcal{I}_t$ non-empty for every $t\in [0, t_f]$. By Lemma~\ref{thm:bangride}, we may assume $\mathcal{I}_{t_f}\neq \emptyset$ under the optimal control sequence $u^\star$. Assuming
    $$
        \mathcal{I}_t\neq \emptyset,\quad \forall t\in[\tau+1,t_f]
    $$
    for some $\tau \in [0,t_f-1]$, we will prove
    $$
    \mathcal{I}_{\tau}\neq \emptyset.
    $$
    Define a control sequence $\hat{u}(\varepsilon)=\lbrace\hat{u}_t\rbrace_{t=0}^{t=t_f}$ by perturbing the optimal solution $u^{\star}$ as in  (\ref{eqn:uhat_nonlin_decreasing}) where $\varepsilon_{\tau}>0$. Let the perturbations $\varepsilon\in\mathbb{R}^{t_f-\tau+1}$ be chosen such that the active constraints in the range $t\in[\tau+1,t_f]$ remain active under the new control sequence $\hat{u}(\varepsilon)$. This happens if condition (\ref{eqn:dy=0}) holds for all $t\in[\tau+1,t_f]$ and $i\in\mathcal{I}_t$, which in the linear case (\ref{eqn:lin}) takes the form
    \begin{align}\label{eqn:dy=0_lin}
    d\bar{y}_{i,t}&=\bar{D}_i\varepsilon_t
    +\sum_{k=\tau}^{t-1} \bar{C}_iA^{t-k-1}B \varepsilon_{k}\nonumber\\
    &=\sum_{k=\tau}^{t} \bar{g}_{i,t-k}\varepsilon_k=0 
    \end{align}
    as $\varepsilon\to 0$. Summation by parts in (\ref{eqn:dy=0_lin}) gives
    \begin{equation}\label{eqn:gsigma}
 \bar{g}_{i,0}\sigma_{t-\tau}=
    \sum_{k=0}^{t-\tau-1} (\bar{g}_{i,t-\tau-k-1}-\bar{g}_{i,t-\tau-k})\sigma_k,
    \end{equation}
    where
    $$
    \sigma_k=\sum_{j=\tau}^{\tau+k}\varepsilon_j.
    $$
    {Since the step response of (\ref{eqn:system}) is concave, the impulse response of the augmented system (\ref{eqn:gt}) is decreasing in the range $t\in[0,t_f]$ ($\bar{g}_{i,t}$ is a decreasing function of $t$). In addition, $\bar{g}_{i,0}=\bar{D}_i>0$ (by Assumption~\ref{ass:h}) and $\sigma_0=\varepsilon_{\tau}> 0$. Therefore,} setting $t=\tau+1$ in (\ref{eqn:gsigma}) gives $\sigma_1 \geq 0$. Similarly, an induction on $t=\tau+2,\tau+3,\dots,t_f$ yields
    \begin{equation}\label{eqn:sigma>0}
            \sigma_k \geq 0, \quad k=0,1,\dots,t_f-\tau
    \end{equation}
    from (\ref{eqn:gsigma}). The first-order perturbation of the cost function is given by (\ref{eqn:Cost_purt}) which by assuming linearity (\ref{eqn:lin}) takes the form
    \begin{align}\label{eqn:Jsigma} 
    \Delta(\epsilon)&=J(x_0,\hat{u}(\varepsilon))-J(x_0,u^{\star})=\nonumber\\
    &=\sum_{t=\tau}^{t_f} F\varepsilon_t+\sum_{k=0}^{t_f-\tau-1}\sum_{t=\tau}^{t_f-k-1}EA^k
    B\varepsilon_{t}\nonumber\\
    &=F\sigma_{t_f-\tau}+\sum_{k=0}^{t_f-\tau-1}EA^k
    B\sigma_{t_f-k-1-\tau},
    \end{align}
    where all the matrices are non-negative by Assumptions~\ref{ass:L} and \ref{ass:f}. As all $\sigma_k$'s are non-negative, the summands in (\ref{eqn:Jsigma}) are non-negative and 
    \begin{equation}\label{eqn:Deltalin>=0}
    \Delta(\varepsilon)\geq 0
    \end{equation}
    holds true. The rest of the proof is identical to that of Theorem~\ref{thm:nonlin_decreasing} and hence omitted.
\end{proof}



\section{HEURISTIC SELECTOR CONTROL STRATEGIES} \label{sec:selector_pid}

As discussed in the previous sections, although closed-loop policies are more robust to model uncertainties than open-loop policies, they require full state measurements to generate the control signal, which can limit their value in applications where only certain outputs are measured by sensors. The bang-ride control policy (\ref{eqn:selector}), in particular, has an interesting structure that can partially address this problem. This policy uses $p+1$ controllers, each designed to keep a corresponding signal at a constraint boundary, and activates (selects) the controller that produces the smallest control signal. In an output feedback setting, where only the process outputs are measured, it is logical to retain this structure. When we have a good model of the process dynamics, we can, for example, estimate the state vector using an observer and replace the true states in (\ref{eqn:selector}) with estimated ones. When we lack a good process model, on the other hand, we may try to replace the state feedback laws in (\ref{eqn:selector}) with linear output feedback controllers, such as PID controllers, leading to a selector control system.

\begin{figure}
	\begin{center}
    \includegraphics[width=1\linewidth]{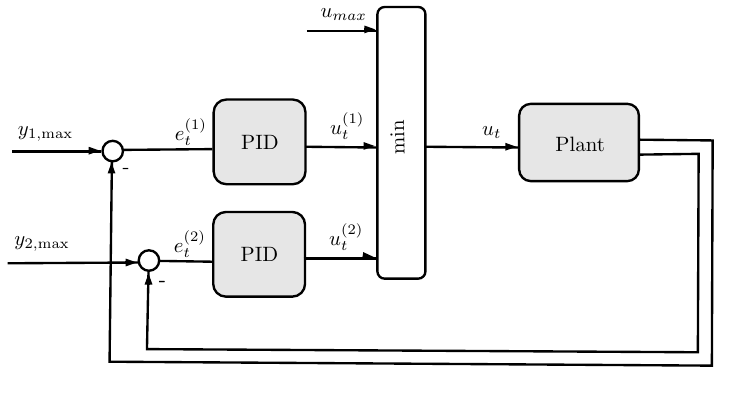}  
	\caption{Selector control structure for a system with $p=2$. The anti-windup feedback links are not shown in this figure.}
	\label{fig:selector}
	\end{center}
\end{figure}

Figure~\ref{fig:selector} shows such a selector control system for (\ref{eqn:system}) which has one input and $p=2$ outputs. There are a total of $p$ PID control loops, one for each output in (\ref{eqn:system}). However, only one loop is connected at each time step. The active loop is chosen by the min selector device. This separation of output channels simplifies the controller design significantly, as one can tune each set of PID parameters separately based on the isolated single-input-single-output feedback loop.

There are several considerations for tuning the PID controllers in a selector control system that approximates a bang-ride control policy. Among the available PID tuning methods, the rules that are based on observing a system response are particularly attractive, as making a single step in the input signal and observing the $p$ outputs allows tuning all controllers based on only one experiment. For the design specification, special attention is paid to setpoint following. The reason is that the occurrence of an overshoot in the $i^{\rm th}$ loop implies that $y_{i,t} > y_{i,{\rm max}}$, \emph{i.e.}, a constraint violation in (\ref{eqn:optimization}). Moreover, a slow response is also undesirable as it leads to large optimality gaps $\vert y_{i,{\rm max}}-y_{i,t}\vert$ ($i\in\mathcal{I}^{br}_t$) when the bang-ride policy is optimal. Based on the above arguments, a suitable control design technique for the selector control loops should minimize overshoots and seek a fast response. Several tuning methods are available that account for these considerations~\cite{pidsurvey}. One such method is presented in~\cite{mtns}, for example, which relies on the first few impulse-response samples of the plant to provide a set of gains that promote the positivity of the tracking errors $y_{i,\rm max}-y_{i,t}$ and reduce overshoot.


Another complication that arises in selector control systems is integrator windup in the PID blocks when the loops switch. However, this problem can be easily treated using a back-calculation technique~\cite[\S3.5]{pidbook}, which results in the following control signals generated by the PID controllers
\begin{align*}
u^{(i)}_t&=k^p_{i}e^{(i)}_t
+k^I_{i} \sum_{k=0}^{t-1}e^{(i)}_k
+k^d_i(e^{(i)}_t-e^{(i)}_{t-1})\\
&+k^w_{i} \sum_{k=0}^{t-1} u_k-u^{(i)}_k
\end{align*}
for all $i=1,2,\dots,p$, where $e^{(i)}_t=y_{i,\rm max}-y_{i,t}$ and $k^w_i$ is the anti-windup gain. When $k^w_i=0$ there is no anti-windup in operation and the $i^{\rm th}$ controller falls back to a basic PID controller. The following control input $u$ is fed to the system (\ref{eqn:system}) from the selector device
\begin{equation}\label{eqn:umin}
    u_t=\min\lbrace u_{\rm max},u^{(1)}_t,\dots,u^{(p)}_t\rbrace,\quad t\in\mathbb{N}_0.
\end{equation} 

\section{APPLICATION TO FAST CHARGING}\label{sec:applications}
In this section, we will demonstrate how the results developed in the previous sections can be applied to the design of fast-charging protocols for lithium-ion batteries in a way that limits their degradation over time. This problem has received significant attention in the research community recently; see the review~\cite{fastcharging}.

The constant-current-constant-voltage (CCCV) protocol typically used for charging lithium-ion batteries in the industry is of bang-ride form: It starts charging the battery with the maximum current and then adapts the charging current to keep the voltage across the battery constant. The bang-ride behavior is also observed in many optimization-based protocols derived for fast charging batteries with different models. For example, a bang-ride solution emerges from the application of nonlinear model predictive control to a fast-charging problem that involves a multi-physics circuit battery model in~\cite{LVWWZTWZ2021}. Similarly, bang-ride solutions were proven to be optimal for three different single-particle battery models in~\cite{PLATM2020} as well as for an equivalent-circuit battery model in~\cite{IFAC}. More generally, it has been shown that under appropriate conditions, bang-ride solutions satisfy necessary conditions of optimality in a wide range of fast charging problems~\cite{shi2025bang}. 

Despite the available results that suggest bang-ride solutions as optimal fast charging protocols, these solutions are not optimal for every optimization problem that describe fast charging lithium-ion batteries. That depends on the type of objectives, constraints, and battery models. This has been documented in~\cite{klein2011optimal,shi2025bang} for example.

Nevertheless, it has been noted that most optimization problems used to obtain health-aware fast charging protocols for lithium-ion batteries satisfy monotonicity assumptions similar to those stated in Section~\ref{sec:monotone_optimal_control_problems} and can therefore be treated as monotone optimal control problems~\cite{DCHCG2023}. Hence, to determine whether a fast-charging problem has an optimal solution in bang-ride form or not, one can use the results of Section~\ref{sec:BR}. When bang-ride solutions are found to be optimal, they can be implemented using an exact state feedback law or an approximate selector control system. We demonstrate this approach on two different battery models with general parameters in this section. Note that the available results in the literature prove the optimality of bang-ride protocols for certain fast charging problems with specific models and parameters, or provide necessary conditions that suggest the optimality of bang-ride protocols for general fast charging problems. This paper fills the gap in this area by providing sufficient conditions that prove the optimality of bang-ride protocols for general fast charging problems with general model parameters.

\begin{figure}
	\begin{center}
    \includegraphics[width=1\linewidth]{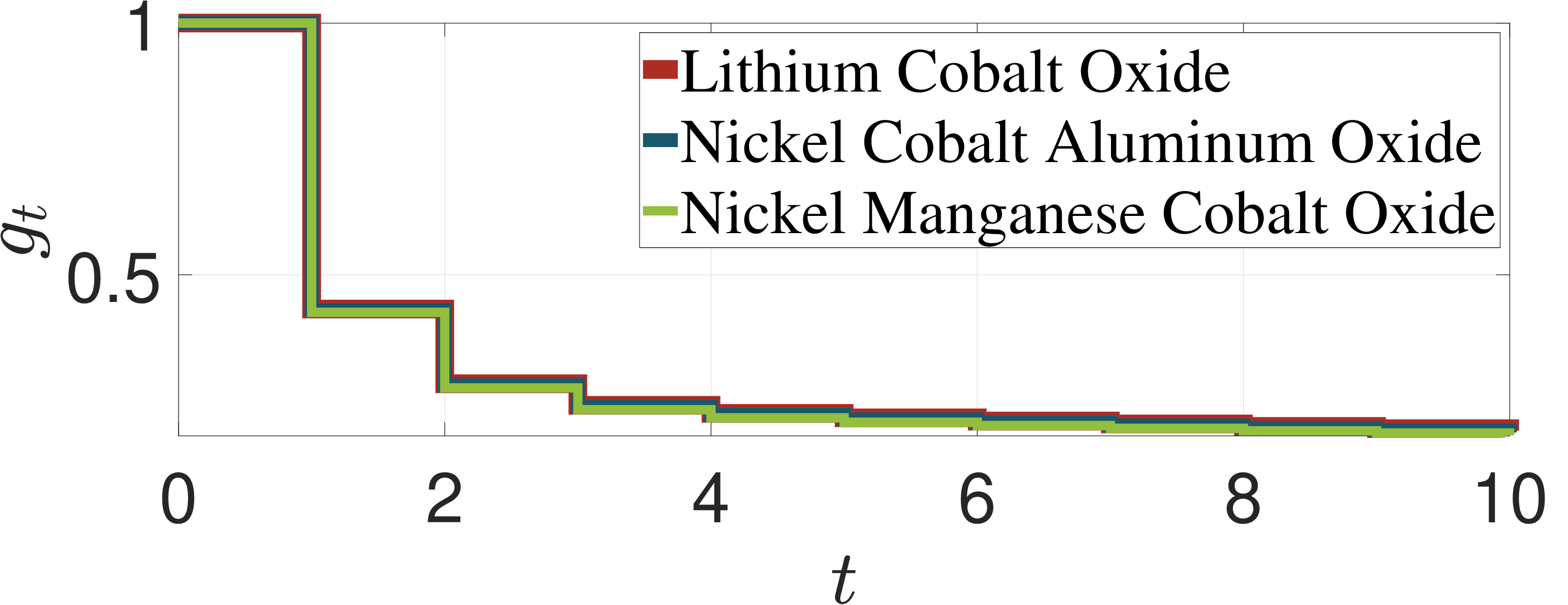}  
	\caption{Normalized impulse responses of three batteries with different chemistries~\cite{PLATM2020,park2}.}
	\label{fig:gt}
	\end{center}
\end{figure}

{\subsection{Single-Particle Models (SPM)}
Single-Particle Models (SPM) are perhaps the most common reduced-order electrochemical models used for lithium-ion batteries. A continuous-time SPM derived from a third-order Pad{\'e} approximation of the diffusion equation that describes the intercalation process and mass transport in the electrodes is given by the following state-space equations~\cite{PLATM2020}
\begin{align}\label{eqn:SPM}
\dot{x}_1(t)&=-a_1x_1(t)+b_1u^-(t) \nonumber\\
\dot{x}_2(t)&=-a_2x_2(t)-b_2u^-(t) \nonumber\\
\dot{x}_3(t)&=-b_3u^-(t) \nonumber\\
y(t)&=-c_1 x_1(t)+c_2 x_2(t)+c_3 x_3(t),
\end{align}
where $u^-$ is the discharging current, $x_3$ is the bulk concentration for the anode, $y$ is the surface concentration, and all the coefficients $a_1$, $a_2$, $b_1$, $b_2$, $b_3$, $c_1$, $c_2$, $c_3$ are positive. 
After the change of coordinates suggested in~\cite{DCHCG2023} and defining the control input as $u(t)=-u^-(t)$, we obtain the model
\begin{align*}
\dot{x}_1(t)&=-a_1x_1(t)+b_1u(t) \nonumber\\
\dot{x}_2(t)&=-a_2x_2(t)+b_2u(t) \nonumber\\
\dot{x}_3(t)&=b_3u(t) \nonumber\\
y(t)&=c_1 x_1(t)+c_2 x_2(t)+c_3 x_3(t).
\end{align*}
Euler discretization results in the discrete-time model
\begin{align}\label{eqn:SPM_DT}
\begin{array}{rcl}
    x_{t+1}&=& Ax_t+Bu_t \\
    h(x_t,u_t)&=& Cx_t 
\end{array}& \quad  \mbox{ for } t\in\mathbb{N}_0,
\end{align}
where
\begin{align*}
    A&=\begin{bmatrix}
    1-a_1t_s& 0& 0 \\ 0 &1-a_2t_s& 0 \\ 0& 0& 1
\end{bmatrix},\;
B=\begin{bmatrix}
    b_1t_s \\ b_2t_s \\ b_3t_s
\end{bmatrix},\;\\
C&=\begin{bmatrix}
    c_1 & c_2 & c_3
\end{bmatrix}.
\end{align*}
A fast-charging problem was proposed in \cite{PLATM2020} based on (\ref{eqn:SPM}) to maximize the bulk concentration at the anode while restricting both the surface concentration and the charging current. For the discretized model (\ref{eqn:SPM_DT}), this optimal control problem can be written as
\begin{equation}
    \begin{array}[c]{rll}
    \underset{u}{\text{maximize}} & J(x_0,u)=\sum_{t=0}^{t_f} x_{3,t} \\
    \mbox{subject to} & y_{t}\leq y_{\max},\\
    & u_t \leq u_{\rm max},\\
    & x_{t+1}=Ax_t+Bu_t, \\
    & y_t=\hat{C}x_t+\hat{D}u_t.
    \end{array} \label{eqn:opt_SPM}
\end{equation}
Here, since the system model (\ref{eqn:SPM_DT}) does not have a direct feedthrough to the outputs, we have used the approach in (\ref{eqn:hhat}) to redefine the output functions as
%
\begin{equation}\label{eqn:y=CAx+CBu}
y_t=\hat{h}(x_t,u_t)=CA x_t + CB u_t:=\hat{C}x_t+\hat{D}u_t. 
\end{equation}
Since $\hat{D}=CB> 0$, Assumption~\ref{ass:h} is satisfied in (\ref{eqn:opt_SPM}). Moreover, by choosing a sampling time in the range
\begin{equation}\label{ts:SPM}
    t_s\in(0,1/\max \lbrace a_1,a_2\rbrace],
\end{equation}
Assumption~\ref{ass:f} is also satisfied. 
Therefore, the optimal control problem~(\ref{eqn:opt_SPM}) satisfies Assumptions~\ref{ass:L}, \ref{ass:f}, \ref{ass:h} and, according to Theorem~\ref{thm:Linear_bangride}, the optimal charging protocol in (\ref{eqn:opt_SPM}) is bang-ride if {the step response is concave, or equivalently, if the impulse response}
\begin{equation}\label{eqn:g_t_hat}
g_t=\left\lbrace\begin{array}{ll}
     \hat{D},&  t=0\\
     \hat{C}A^{t-1}B,& t\geq 1 
\end{array}\right.
\end{equation}
is decreasing in the range $t\in[0,t_f]$. This is true for all the three sets of battery parameters considered in~\cite{PLATM2020,park2}, as shown in Figure~\ref{fig:gt}. Therefore, the optimal charging protocol for these battery models is bang-ride, a result that is consistent with \cite{PLATM2020}.
In contrast with \cite{PLATM2020}, our results allow to prove optimality for \emph{all} SPM battery model parameters $a_1$, $a_2$, $b_1$, $b_2$, $b_3$, $c_1$, $c_2$, $c_3$ that result in a decreasing impulse response $g_t$.}

\subsection{Equivalent-Circuit Models (ECMs)}\label{sec:ECM}

Equivalent-Circuit Models (ECMs) are one of the most popular empirical models of lithium-ion battery cells~\cite{ECM}. A continuous-time ECM with two time constants is given by the following set of ordinary differential equations
\begin{align}\label{eqn:ECM1}
\dot{v}_1(t)&=-\frac{1}{R_1C'}v_1(t)+\frac{1}{C'}u(t) \nonumber\\
\dot{v}_2(t)&=-\frac{1}{R_2C''}v_2(t)+\frac{1}{C''}u(t) \nonumber\\
\dot{s}(t)&=\frac{1}{Q}u(t) \nonumber \\
y(t)&=v_{1}(t) + v_{2}(t) +\beta s(t)+R_0 u(t).
\end{align}
Here, $u$ is the charging current, $v_1$, $v_2$ are the voltages across the (abstract) RC-links, $s$ is the state of charge, and $y$ is the over-potential voltage. This model is based on the assumption that the open-circuit voltage function of the battery is known and the state of charge can be computed via integration of the input current $u$. A forward Euler discretization of \eqref{eqn:ECM1} with sampling time $t_s$ yields the following discrete-time battery model
\begin{align}\label{eqn:ECM_DT}
    x_{1,\,t+1}&=(1-t_s/(R_1C'))x_{1,t} 
    +\frac{t_s}{C'}u_t
    \nonumber\\
    x_{2,\,t+1}&=(1-t_s/(R_2C''))x_{2,t}
    +\frac{t_s}{C''}u_t
    \nonumber\\
    x_{3,\,t+1}&=x_{3,t}
    +\frac{t_s}{Q}u_t
    \nonumber\\
    h(x_t,u_t)&=x_{1,t} + x_{2,t} +\beta x_{3,t}+R_0 u_t.
\end{align}
We assume $x_0=0$ as the initial condition. Now consider the fast charging problem (\ref{eqn:optimization}) with the cost function
\begin{equation}\label{eqn:L=s_t}
    L(x_t,u_t)=x_{3,t}.
\end{equation}
We consider constraints on the cell charging current $u_t$ as well as its over-potential voltage $y_t=h(x_t,u_t)$ to promote safety and limit degradation~\cite{kabir}. In particular, constraining the over-potential rather than the terminal voltage (as in conventional CCCV protocols) helps prevent stress and strain~\cite{IFAC}. These constraints can be written as
\begin{align*}
u_t&\leq u_{\rm max}, \\
y_{t}&\leq y_{\rm max}
\end{align*}
for some positive constants $u_{\rm max}$ and $y_{\rm max}$. By choosing a sampling time $t_s \in (0, t_s^{(1)}]$, where
\begin{equation}\label{eqn:Ts1}
    t_s^{(1)} =\min\lbrace R_1C',R_2C'' \rbrace,
\end{equation}
the transition function $f(x_t,u_t)$ is increasing in both $x_t$ and $u_t$. It is also obvious that $f(x_t,u_t)$ is continuous everywhere. Therefore, Assumption~\ref{ass:f} is satisfied. It can similarly be shown that Assumptions~\ref{ass:L} and \ref{ass:h} are also satisfied. Therefore, the fast charging problem described above is a monotone optimal control problem, and its optimal solution satisfies the input maximization property described in Lemma~\ref{thm:bangride}. By further restricting the sampling time to
$t_s\in(0,\min\lbrace t_s^{(1)},t_s^{(2)}\rbrace]$,
where
\begin{equation}\label{eqn:Ts2}
t_s^{(2)}=\frac{R_0}{\frac{1}{C'}+\frac{1}{C''}+\frac{\beta}{Q}},
\end{equation}
the (more restrictive) hypothesis of Theorem~\ref{thm:Linear_bangride} is also satisfied. To see this, we first note that both the cost function~(\ref{eqn:L=s_t}) and the process dynamics (\ref{eqn:ECM_DT}) are linear with the impulse response satisfying $g_0=R_0$ and
$$
g_t= \frac{\left(1-\frac{t_s}{R_1C'}\right)^{t-1}}{C'/t_s} +
\frac{\left(1-\frac{t_s}{R_2C''}\right)^{t-1}}{C''/t_s}+
\frac{t_s\beta}{Q}
$$
for $t\geq 1$. Inequality $0<t_s\leq t_s^{(2)}$ implies
$$
g_0=R_0\geq \frac{t_s}{C'}+\frac{t_s}{C''}+\frac{t_s\beta}{Q}=g_1,
$$
while inequality $0<t_s\leq t_s^{(1)}$ ensures that the impulse response $g_t$ is decreasing in the range $t\in[1,t_f]$ for any $t_f\geq 1$. Therefore, when both inequalities are met, \emph{i.e.},
$$
t_s\in(0,\min\lbrace t_s^{(1)},t_s^{(2)}\rbrace],
$$
the monotone optimal control problem satisfies the hypothesis of Theorem~\ref{thm:Linear_bangride}. This indicates that the optimal solution is of bang-ride form. Therefore, to implement the exact optimal closed-loop policy, one may solve the constraint equations (\ref{eqn:constraint_equation}) analytically as follows
\begin{align}\label{eqn:Kbangridebattery}
K_0(x_t)&=u_{\rm max},\\
K_1(x_t)&=(y_{\rm max}-x_{1,t}-x_{2,t}-\beta x_{3,t})/R_0.\nonumber
\end{align}
Alternatively, the optimal closed-loop policy can be implemented approximately, using a selector control system based on PID controllers (see Section~\ref{sec:selector_pid} and Figure~\ref{fig:selector}). 

\begin{figure}
     \centering
     \begin{subfigure}[b]{0.48\textwidth}
        \centering
        \includegraphics[width=1\linewidth]{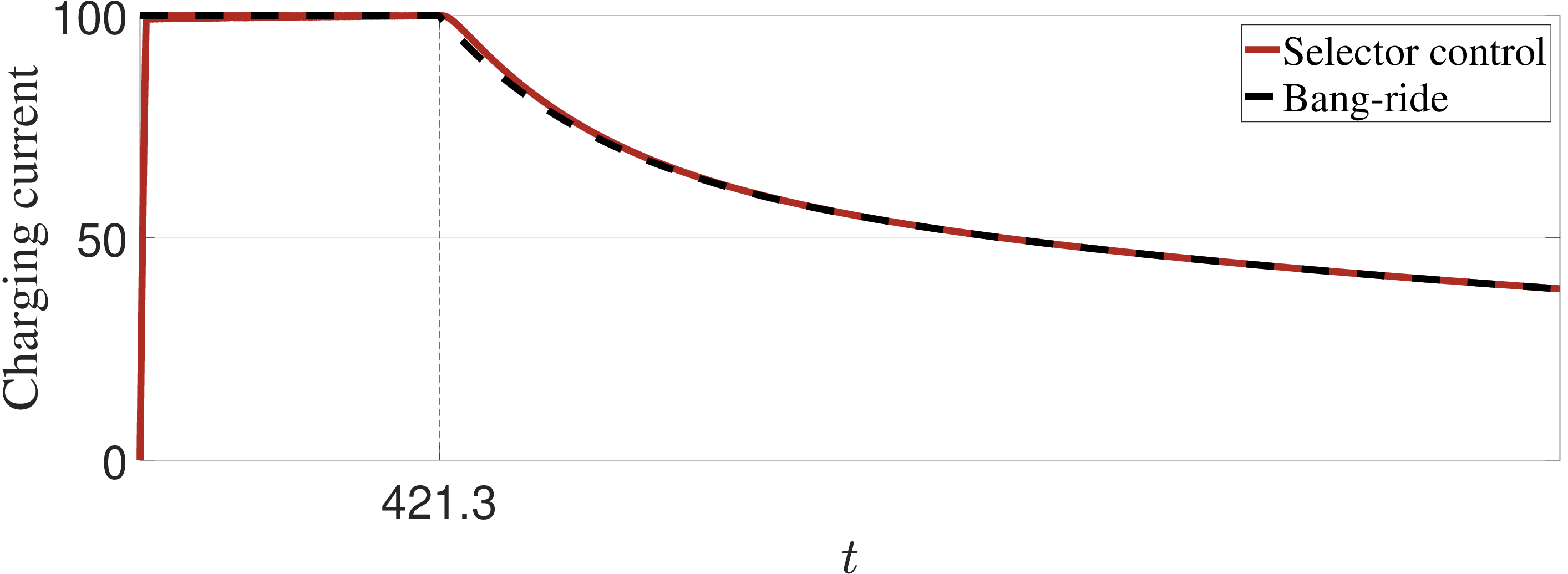}
	    \caption{Charging current over time.}
         \label{fig:ECM_i}
     \end{subfigure}
     \hfill
     \begin{subfigure}[b]{0.48\textwidth}
         \centering
        \includegraphics[width=1\linewidth]{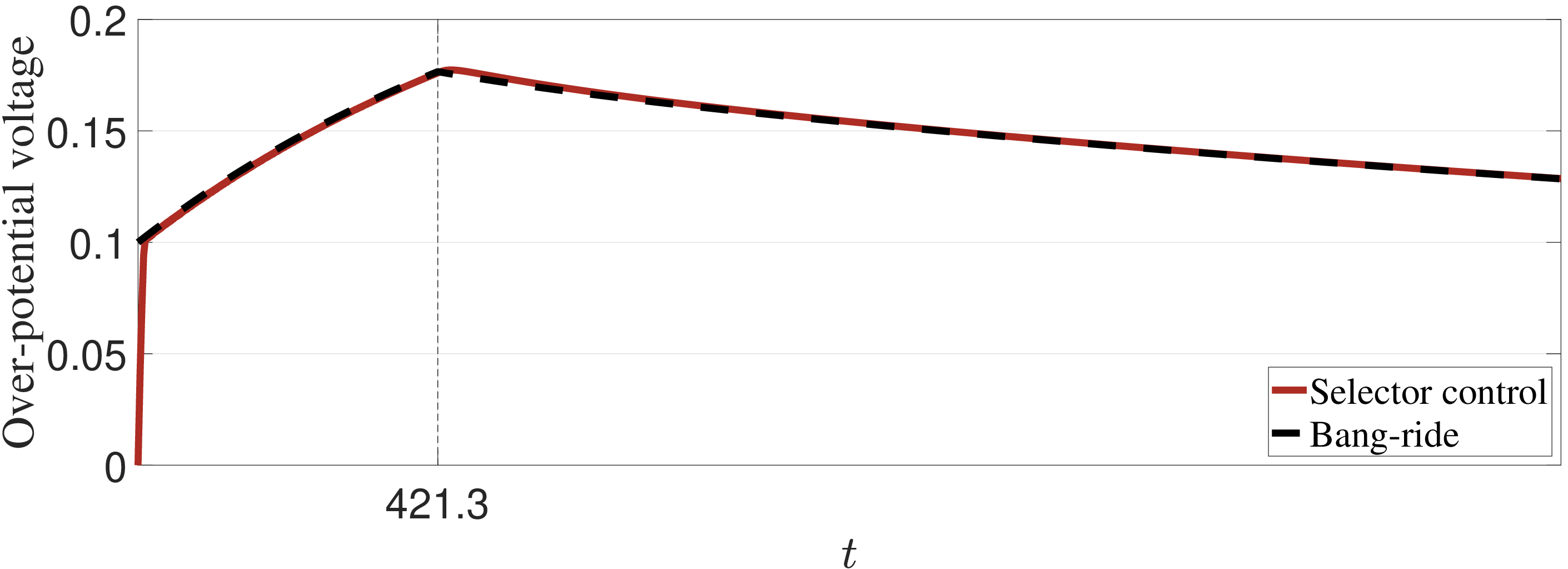}
	    \caption{Over-potential voltage over time.}
         \label{fig:ECM_v}
    \end{subfigure}
    \hfill
    \begin{subfigure}[b]{0.48\textwidth}
         \centering
        \includegraphics[width=1\linewidth]{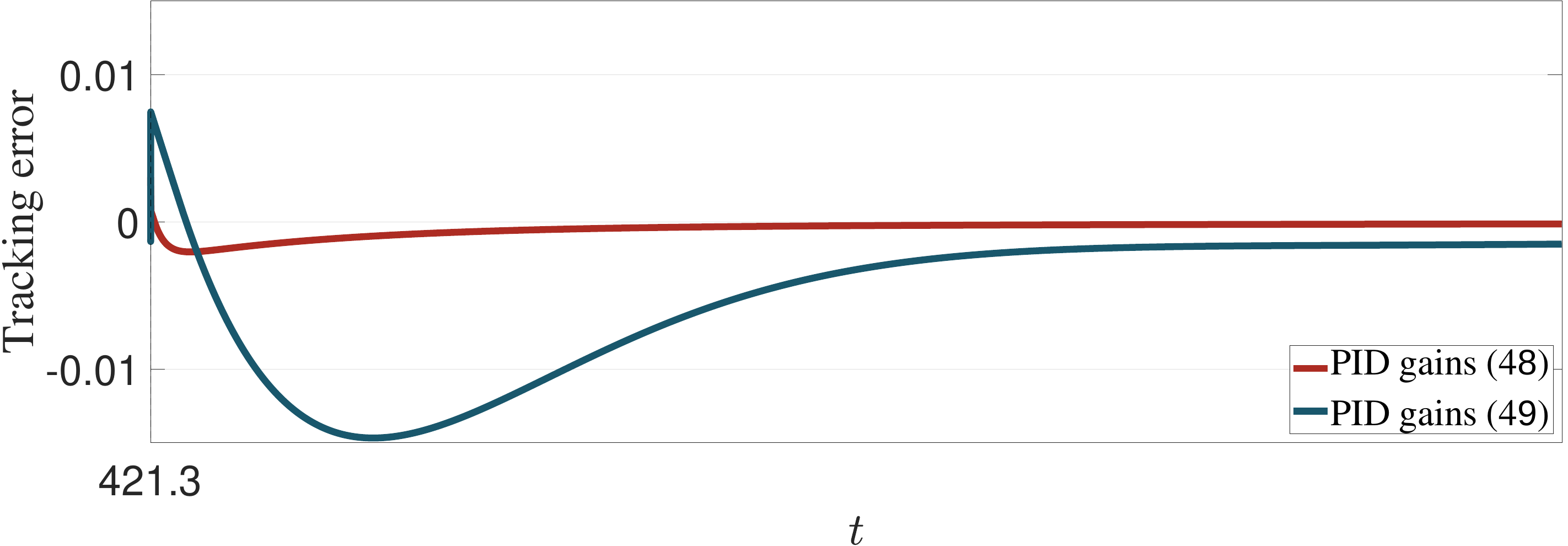}
	    \caption{Active constraint tracking errors over time.}
       \label{fig:ECM_e}
    \end{subfigure}
        \caption{Implementation of the optimal charging profile of the constrained fast charging problem considered in Section~\ref{sec:ECM} under two different strategies: exact implementation via the state-feedback control law (\ref{eqn:selector}) and an approximate implementation via a selector control system shown in Figure~\ref{fig:selector}. 
        }
        \label{fig:ECM1}
\end{figure}

To simulate the optimal charging protocol we consider the following parameters for the battery model
$$
\left\lbrace\begin{array}{l}
R_0=1 \,(m\Omega),\; R_1=1.5\, (m\Omega),\; R_2=1\, (m\Omega)\\
C'=2000 \,(kF),\; C''=500\, (kF),\; Q=50 \,(Ah)\\
\beta=0.1,\;u_{\rm max}=100 \,(A),\; y_{\rm max}=0.2\, (V) 
\end{array}\right..
$$
The sampling time limits (\ref{eqn:Ts1}) and (\ref{eqn:Ts2}) are given by $t_s^{(1)}=500\, (s)$ and $t_s^{(2)}=327.27\, (s)$, respectively, which are well above the chosen sampling time $t_s=0.05 \,(s)$. The resulting bang-ride charging current in (\ref{eqn:Kbangridebattery}) and the charging current $u$ of the heuristic selector control scheme (\ref{eqn:umin}) are shown in Figure~\ref{fig:ECM_i}, as dashed and solid plots respectively. In Figure~\ref{fig:ECM_v}, the over-potential voltages of the cell are also compared for these two charging profiles. The selector control scheme relies on the measured outputs and PID controllers, whereas the bang-ride charging current (\ref{eqn:selector}) is based on the exact battery model and requires accessing the state. Nevertheless, the two charging profiles have almost indistinguishable performances in this experiment, as seen in Figure~\ref{fig:ECM1}. The main deviation from the optimal solution occurs when the active constraint is switched from current to voltage at $t=421.3$. In this experiment, we have used the design technique from~\cite{mtns} to find the PID gains
\begin{align}\label{pidgains_mtns}
    k^p_1=1,\, k^i_1=k^w_1=4,\, k^d_1=0.5.
\end{align}
This set of gains was compared with the gains
\begin{align}\label{pidgains_matlab}
    k^p_1=43.23,\, k^i_1=k^w_1=0.366,\, k^d_1=0
\end{align}
suggested by MATLAB built-in tuner. Figure~(\ref{fig:ECM_e}) shows the errors that occur after constraint switching under the two different controller tunings. We observe that the gains (\ref{pidgains_mtns}) result in a significantly smaller constraint violation. Whether or not there is a combination of PID gains and anti-windup procedures that guarantees feasibility at all times is not known in general.

\section{CONCLUSIONS}\label{sec:conclusion}
This paper considered monotone optimal control problems governed by monotone costs, constraints, and dynamics. Such problems appear in many applications, including in the design of fast-charging protocols used for lithium-ion batteries. Optimality conditions were obtained for several classes of monotone optimal control problems, and families of systems that admit bang-ride optimal solutions were identified. It was shown that the corresponding optimal control policies can be implemented by a selector control strategy that switches between a finite number of state feedback laws. These structural results motivate the use of selector strategies based on PID controllers that are typically designed without perfect models and full-state measurements. Simulations of single-particle models and equivalent-circuit models of lithium-ion batteries in fast-charging problems demonstrated the efficacy of both optimal and heuristic selector control strategies.

\bibliography{ref.bib}
\bibliographystyle{ieeetr}

\end{document}